
\def\g{\gamma}     \def\l{\lambda}    \def\o{\omega}   
\def\D{{I\!\!D}}     \def\p{\partial}
\def\e{\varepsilon}      \def\C{{\bar C\!\!\!\!I}} 
\def\cl{{\rm {cl}}}   \def\grad{{\rm {grad}}}   \def\id{{\rm {id}}}
\def\Arg{{\rm {Arg}}}   \def\Comp{{\rm{Comp}}}   \def\Teich{{\rm{Teich}}}
\def\Crit{{\rm{Crit}}}

\centerline{\bf  Iterations of rational functions: }    

\centerline{{\bf which hyperbolic 
components contain polynomials?}\footnote{$^1$}{This 
is a revised and 
extended version of a 
 manuscript "Polynomials in hyperbolic components" 
 IMS SUNY at Stony 
Brook, June 1992}}

\


\

\centerline{ by Feliks Przytycki\footnote{$^2$}{Research partially supported by Polish KBN Grants 210469101
"Iteracje i Fraktale" and 210909101 "...Uklady Dynamiczne". }}

\

\

{\bf Abstract.}  \  {\it Let  $H^d$ be the set of all rational 
maps of degree $d\ge 2$ on the Riemann sphere
which are expanding on Julia set. We prove that if $f\in H^d$ 
and all or all but one critical points (or values) are in the immediate
basin  of attraction to an attracting fixed point 
then there exists a polynomial in the component $H(f)$ 
of  $H^d$ containing $f$. If all critical points are in the immediate basin 
of attraction to an attracting fixed point or 
parabolic fixed point then $f$ restricted to Julia set is conjugate  
to the shift on the one-sided shift space of $d$ symbols. 

We give {\rm exotic} examples of maps of an 
arbitrary degree $d$
with a non-simply connected,  completely invariant 
basin of attraction 
and arbitrary number $k\ge  2$ of critical points in the basin. For such a
map $f\in H^d$ with $k<d$ there is no polynomial in $H(f)$. 

Finally we describe a computer experiment joining an exotic example to 
a Newton's method (for a polynomial) 
rational function with a 1-parameter family of rational 
maps.}

\

\

{\bf Introduction.}

\

In the space $Q^d$ of rational maps of degree $d\ge 2$ of the 
Riemann sphere $\C$, denote by $H^d$ the set of maps which are expanding on
the Julia set. Expanding means that there exists $n>0$ such that for every $z$ 
in the Julia set 
$|(f^n)'(z)|>1$ (the derivative is in the standard spherical 
Riemann metric). We call $z\in \C$ a critical  point if 
$f'(0)=0$. We call $v$ a critical value if $v=f(z)$ for 
a critical point $z$.

In Section 1 we prove following:

\

{\bf  Theorem A.} Let $f\in H^d$. Suppose that all, 
or all but one, of the critical values of $f$ are 
in an immediate basin of attraction $B(f)$ to one attracting $f$-fixed point
$p$.  Then the component  $H(f)$ of $H^d$  containing $f$ contains also a 
polynomial.

\

(Critical values are counted in Theorem A without multiplicities. However
critical  
points everywhere in the paper are counted with multiplicities.)

\

{\bf Corollary B.} If all critical points of $f$ are in 
$B(f)$ then $f$
restricted to the Julia set $J(f)$ is conjugate to the full one-sided  shift.

\

{\bf Theorem C.}  Let $f\in Q^d$, \  $p$ be a parabolic 
fixed point, and let $B(f)$ be an immediate basin of attraction 
to $p$ adjacent to $p$ such that $f(B(f))=B(f)$ (which is equivalent to $f'(p)=1$). Suppose that all critical points of 
$f$ are in $B(f)$. Then 
$f$ 
restricted to  $J(f)$, is conjugate to the full one-sided  shift.

\

The procedure to prove Theorem C is similar to 
Corollary B so it will be only sketched. 

\

It is much easier to prove that,  under the assumptions of 
Corollary B or 
Theorem C, $J(f)$ is Cantor, $f|_{J(f)}$ is conjugate to a topological 
Markov chain and $f^n|_{J(f)}$ is conjugate to a 1-sided shift, 
then to prove that $f|_{J(f)}$ itself is conjugate 
to a full 1-sided shift.

\

The questions answered in Corollary B and Theorem C  were asked to me by John Milnor. In the case of the basin of a sink 
he suggested  to join critical values with the sink along trajectories of 
the gradient flow , Morse curves, described below.  
This was a fruitful idea. After proving Theorem A 
and Corollary B I learned that these facts were proved by  
P. Makienko\footnote{$^3$}{Recently his proof appeared in 
the preprint  [M]. It influenced our revised version of the paper. Also our Theorem C is proved in [M].} in his PhD paper but stayed unpublished. Corollary B is 
also stated in [GK] but proved only for $d=2$.

\

Pay attention that Corollary B and Theorem C depend 
on holomorphic 
phenomena. Indeed there exists a 1-sided 
topological Markov chain $T$ 
for which each point has 10 preimages which has the same 
$\zeta$-function as the full 1-sided shift of 10 symbols $S_{10}$ 
and which sufficiently high power $T^n$ is conjugate to $S_{10}^n$, 
but $T$ is not conjugate to $S_{10}$ [Boyle].

\

If $B$, the immediate basin of attraction to an 
attracting  fixed point,
 is simply connected, then the 
number of critical 
points of $f$ in $B$ is equal to $\deg (f|_B)-1$ . 
(Because $f$ pulled-back to the unit disc $\D$ by a Riemann 
mapping  is a Blaschke product which has 
$\deg (f|_B)-1$ critical points in $\D$ (and the same 
number of them outside).)

If $B$ is the basin of attraction to $\infty$ for a polynomial 
and $B$ is not simply-connected, then the
number of critical 
points of $f$ in $B$ is at least $\deg (f|_B)$ (including $\infty$ as $(\deg (f)-1)$-multiple critical point ).  

Surprisingly this is false in general.  The "proof" that if the basin is not simply-connected than it contains at 
least as many critical points  as the degree of $f$ on it , given in [P], 
is wrong for degree larger than 2, the corresponding 
Lemma in [P] is false.

Here in Section 2 we prove with the use of the 
quasiconformal surgery technique [D] the existence of 
{\it exotic basins}: 

{\bf Theorem D.} There exists a rational function $f$ 
of an arbitrary degree $d\ge 3$ with a completely 
invariant (i.e. invariant under $f^{-1}$) non-simply 
connected basin of  immediate attraction  to 
an attracting fixed point, and with an arbitrary $2\le k \le 2d-2$   
number of  critical points in the basin. 

\

The conclusion is that the assumption all or all but one 
critical values are in the basin in Theorem A, is essential.  
Namely, due to Theorem D for $k\le d-1$, we arrive at 

{\bf Corollary E.} For every $d\ge 3$ there exists $f\in 
H^d$  with a completely invariant basin of immediate 
attraction to 
an attracting fixed point, such that $H(f)$ contains no 
polynomial.

\

Neither in Theorem A nor in Corollary  E does it matter that 
$f\in H^d$. Only the basin $B(f)$ matters. The  reader 
will find appropriate precise assertions in Sections 1 and 
2.

\

 I checked with the help of computer that  $f(z)=z^2+c+b/(z-a)$ for $c=-3.121092, a=1.719727, 
 b=.3142117$ is an exotic example for $d=3$ which 
existence is asserted in Corollary E. I am grateful to Ben Bielefeld and Scott Sutherland for the help in producing 
 several computer pictures of Julia sets for such $f$'s and related pictures 
 in the parameter space. 
Ben invented the parametrization in which the pictures 
were done. Our 1-parameter families join exotic examples 
of the above type of $d=3$, having two superattracting 
fixed points 
  and a critical point of period 2, 
with functions having three superattracting fixed points. It is easy 
to see (and well-known) that the latter functions must be Newton's 
method 
rational functions for degree 3 polynomials. 

(If $f\in Q^d$ has $d$ superattracting fixed points then in 
appropriate holomorphic coordinates on $\C$ it is 
Newton's for a degree $d$ polynomial. Hint to a proof: 
Change first the coordinates on $\C$ by a homography so that unique repelling 
fixed point becomes $\infty$.)
 
\

\

{\bf Section 1.  Rearranging  critical values. Proofs of 
Theorem A and Corollary B.}

\

Let $B(f)$ be the immediate basin of attraction to 
an attracting fixed point $p$  for a rational map $f\in 
H^d$.
Suppose $f'(p)\not=0$. 

\

We shall make in this Section 
 the following types of perturbations of such maps in $H^d$:

\

{\bf 1. A perturbation along a curve $\g$.} We have in mind here the following 
construction: Suppose there is a  curve $\g=\g(t), t\in [0,1]$ embedded in 
the basin $B(f)$ with $p\notin\g$.  Take a small 
neighbourhood $U$ of $\g, U\subset B(f)$ disjoint 
with a neighbourhood of $p$.
 Let $g_t $ for every $t\in [0,1]$ be a diffeomorphism of 
$\C$ so that $g_t(\g(0))=\g(t)$ and $g_t$ be the identity outside $U$ , 
$g_0=\id$ and $g_t$ smoothly depend on $t$. We obtain the homotopy 
$h_t=g_t\circ f$. Pay  attention that though we called our perturbation 
"along $\g$" we change the map in a neighbourhood of $f^{-1}(\g)$. 

If the following assumption holds:

$$p\notin U ,\ \  h_t^n(U)\to p \leqno (1)$$
then of course the basin of attraction to $p$ is the same for $h_t$ as for $f$.

Then  we construct an invariant measurable $L^\infty$ conformal structure for each $h_t$ as 
follows: We take  the   standard structure  on a small neighbourhood of $p$ 
then we pull-back it  by $h_t^{-n}$ . On the complement  we take the standard 
structure. Now we integrate this structure (we refer to 
Measurable Riemann Mapping Theorem [B], [AB]) and $h_t$ in the new coordinates  
gives  a homotopy $f_t$ through maps in $H^d$. See [D] 
for this technique.

\

{\bf 2. A small $C^1$-perturbation.} If a map $g$ is 
$C^1$ close to $f$ on $U$ such that $\cl U\subset 
B(f)\setminus \{p\}$, $f=g$ outside $U$ and in neighbourhoods of critical points $g$ differs from $f$ only by affine parts, then clearly it is homotopic to $f$ 
through also small perturbations satisfying the same conditions as $g$. The condition (1) holds 
automatically. As before we introduce new conformal structures, integrate them and obtain a perturbation 
homotopic through maps in $H(f)$ to $f$.

\

{\bf 3.  Blaschke type perturbation.}  Let  $U\subset 
B(f)$ be an 
open topological disc containing $p$,  with smooth boundary not containing critical points,  
such that $f(\cl U)\subset U$ 
\ \ $f:U\to f(U)$ is a proper map and 
$d':=\deg f|_U \ge 2$. Then we construct a 1-parameter 
family of maps joining in $H(f)$ the map $f$ to a map 
having $p$ as a $d'$-multiple fixed point as follows: 

Let $R_1, R_2$ be Riemann maps from $U$, respect. 
$f(U)$,  to the unit disc $\D$ such that $R_i(p)=0, i=1,2.$ 
Let $a_1=0,a_2,...,a_{d'}$ be $R_1$-images of $f|_U$-preimages of $p$. Let $B_t=\lambda z
\prod_{i=2}^{d'} {z-ta_i \over 1-t\bar{a_i}z}, |\lambda|=1,  
 t\in [0,1]$. We set $h_t=R_2^{-1}\circ B_t\circ  R_1$. 
Here $\lambda$ is chosen so that $h_1=f$.  It is useful 
to write $h_t=R_1^{-1}\circ g_t \circ R_1$, where 
$g_t=R_1\circ R_2^{-1}\circ B_t$. Change $h_t$ in $U$ 
close to $\partial U$ by a smooth isotopy so that $h_t$ 
and  $f$ coincide on $\partial U$ for all $t$. We extend 
$h_t$ outside $U$ by $f$ to the whole $\C$. As in the previous cases we pull-back the standard conformal 
structure  from $R_1(f(U))$ by $R_1^{-1}$ to $f(U)$,  
extend it by $h_t^{-n}$,  complete on $\C\setminus 
\bigcup h_t^{-n}(f(U))$ with the standard structure and integrate.

\

\

Let  $\Phi_f$ conjugates $f$ to $z\to \lambda z$ where 
$\lambda:= f'(p)\not=0$, in a neighbourhood of $p$ 
(i.e. $\Phi_f f(z)=\lambda\Phi_f(z)$. Extend $\Phi_f$ to 
$B(f)$ by $\Phi_f(z)=\l^{-n}\Phi_f f^n(z)$. Define 
$G_f(z)= |\Phi_f(z)|^2$.

\

If $f\in H^d$ and $g\in H(f)$ then  
we write $p_g$ for the point $z:g(z)=z$ such that $(g,z)$ 
belongs to the  component  of 
the Cartesian product $ H(f)\times\C $ containing 
$(f,p_f)$. There exist $B(g)$ and $\Phi_g$ (provided 
$g'(p_g)\not=0$) as above. 

Now we can formulate 

{\bf Main Lemma.}  For every $f\in H^d$ there exists 
$g\in H(f)$ such that $g'(p_g)\not=0$ and there exists 
$a>0$ such that all critical values of critical points in $B(g)$ are in a component $\partial$ 
of $\{G_g=a\}$ which is 
a topological circle 
separating $p_g$ from Julia set $J(g)$.

\

{\bf Proof.}  By small perturbations (types 2 and 3) 
we assure that the sink $p_f$ is neither a critical 
 point  nor a critical value  for a 
critical point in $B=B(f)$  all critical points in $B$ are simple and their 
 forward trajectories are pairwise disjoint.

At the end it may occur useful also to have ${1\over \pi} 
\Arg f'(p)$ irrational. We assure this by type 3 
perturbation where $\Arg\l$ is the parameter.

\

The critical points of $G=G_f$ are

1)  the fixed point $p$ and its iterated pre-images (these are  minimum 
points with $G=0$),

2)  the critical points of $f$ and their iterated pre-images (these are 
saddle points for  $G$).

Denote the set of points in 1) by $M$ and the set of points in 2) by  $S$.

\

For every $q\in M$ let $A(q)$ denote the basin of attraction to $q$ for the 
flow of  the vector field    $-\grad G$.

 Denote by $r(q)$ the least non-negative integer such that $f^{r(q)}(q)=p$. 
 For every  $z\in A(q)$ or $X\subset A(q)$ write $r(z):=r(q)$ and $r(X):=r(q)$

\

Observe that for every  $z\in B(f)$  there exists a curve $\g$ joining $z$ with 
$p$ consisting of critical points of $G$ and of trajectories of  $ \grad  G$ 
where this field is non-zero (i.e. $\g$ goes from  $z$ to a critical point,, 
say a minimum, then to a saddle, then to a minimum etc. until it reaches $p$. 
Let  $\g(z)$ denote  a curve as above intersecting the minimal possible 
number of $A(q)$'s. Denote the number of  these $A(q)$'s by $s(z)$ or $s_f(z)$ and call 
the curve a {\it Morse curve}.

\

{\bf Observations} 

1.  $r(f(z))=r(z)-1$ if $r(z)\ge  1$

2.  $ s(f(z))\le s(z)$.

\

The observation 2. follows from the fact that $f$ maps trajectories of $ \grad  G$ to trajectories 
of $ \grad  G$.

\

 The plan is now to move all critical 
 values to the same level   $ G=a$  in $A(p)$.

\

We shall do it for each critical point separately so that  we do not  move 
the critical values moved before to the level $a$ in $A(p)$.  We move each critical  
value $f(c)$ step by step so that after each step $s(f(c))$ decreases and  
$f(c)$ is in some $A(q)$. When $f(c)\in A(p)$  we move $f(c)$ along the 
trajectory of $ \grad  G$ to the level $a$ as described in 1.

\

Take $c$ a critical point for $f$.  By a small 
perturbation in a neighbourhood of $c$ we obtain $f(c)\in A(q)$.  We can assume
that $q\not= p$. 

This is correct because

a) The perturbation is above the $G$ level of $f(c)$ so it does not  change a part of the stable manifolds of a saddle to which $f(z) $ might belonged, in a neighbourhood of $f(c)$ and between it and the saddle).

b) The change of coordinates integrating the new invariant conformal structure is conformal there. This is so because  the structure is non-standard only in some places above $a$. Thus the change of coordinates maps the gradient lines of the old $ G$ to the gradient lines 
of the new one and  stable manifolds (separatrices) to separatrices.

\

Now let $\g$ be the part of a Morse curve $\g(f(c))$ joining $f(c)$ with the first saddle $\o\in \g(f(c))$.

Observe that for every $j>0$ we have $f^j(\g)\cap\g=\emptyset$.

Indeed, we have $f^j(\g\cap A(q))\cap\g=\emptyset$ by Observation 1. . By the same observation $f^j(\o)\in  \cl  A(f^j(q))$ cannot  be in $\g\cap A(q)$.  Finally $f^j(\o)\to p$ so $f^j(\o)\not=\o$.

As $f^j(\g)\to p$ we can find $U$ a neighbourhood of  $\g$ so that  $f^j(U)\cap U=\emptyset$ for all $j>0$. Take in $U$ a curve $\g'$ joining $f(c)$ with a point  $z$ close to $\o$, $z\in A(q')$, $s(q')=s(q)-1$. Take care about $\g'$ not containing critical values under iterates $f^j$ of other critical points. 

Now perturb $f$ along $\g'$ as described in 1. Do it with $g_t$ different from identity only in a neighbourhood of $\g'$ small enough not to contain critical values  of other critical points, under iterates $f^j$ . 
The condition (1) is of course satisfied. 

Now let us explain why the new map $f_1$ has the property 

$$
s_{f_1}(f_1(c_{f_1})=s_f  (f(c))-1.           \leqno (2)
$$
 Here $c_{f_1}$ is the old $c$ in new coordinates, it is a critical point for $f_1$. We use  the fact that a part of the domain we changed $f$ to $h_t$ is  the basins $A(f^{-1}(q))$  where by  Observation 2. $s\le s(q)$ . So we did not change $ G$ in the basins, the part of the  Morse curve $\g$ beyond $\o$ goes through, where $s<s(q)$. We changed $f$ also in a neighbourhood of  $f^{-1}(\o)$. This does not change $ G$ below a neighbourhood of $\o$ because $ G (f^{-j}(\o)> G(\o), j=1,2,...$. (Compare this with the argument (a).) This does not change  $G$ below the part of $\g$  beyond a neighbourhood of $\o$ neither, because this change  is close to the set $s\ge s(q)$.

 The change of $h_t$ to $f_t$ does not hurt (2).
Indeed also by the above arguments the new measurable conformal structure is the standard one below $\o$ and the part of $\g$ beyond $\o$. So the change of coordinates maps the gradient lines of the old $ G$ to the gradient lines 
of the new one as in
 (b).

\

When $f(c)\in A(p)$ we first make a small perturbation so that the gradient line of $G$ passing through $f(c)$ 
does not intersect forward trajectories of other critical 
values which are already in $A(p)$. Next we move 
$f(c)$ along $\g$ which is a piece of the gradient line of 
$G$ 
passing through $f(c)$ joining it with $\{G=a\}$ in $A(p)$. 
We succeed because ${1\over\pi}\Arg (f^n)'(p)$ is 
irrational so all the curves $f^n(\g), n\ge 0$ are pairwise 
disjoint and the conformal structure does not change 
below these curves.

\ 

After a sequence of consecutive  perturbations as above  we obtain a rational mapping $g$ with all the  critical values  on one  level $\{ G=a\}$, more precisely its 
component  $\partial$ 
intersecting $A(p)$.

Now denote the domain of $\C\setminus \partial$ containing $p$ 
by $D_a$. 
To finish Proof of Main Lemma take 
$b>0$ so close to $0$ that a component 
$\partial '$ of $\{G=b\}$ is in the domain around 
$p$ where $g$ is linearizable. So $\partial '$ is a 
topological circle. Then each $x\in \partial '$ 
can be mapped to the point of intersection of  $\partial $ 
with the trajectory of $\grad G$ starting from $x$. 
This gives a homeomorphism between 
 $\partial '$ and $\partial $.

Otherwise  $\cl D_a$ contains an $S$-type 
 critical point $x$ for  $G_g$. Then there exists $n>0$ such that  $y=g^n(x)$ is a critical value for $g$. Hence 
$G_g(y)=|g'(p_g)|^n  G_g(x)< G_g(x)\le a$. This 
contradicts just  achieved $G_g(y)=a$.

 Proof of Main Lemma is finished.  \hfill$\clubsuit$

\

{\bf Proof of Theorem A.} Let $f$ be  already as $g$ 
 in Main Lemma. 
By perturbing along curves one obtains additionally 
all critical values of critical points from $f^{-1}(B(f))
\setminus B(f)$ also in $\partial$. ({\it A posteriori} 
we will see that under the assumptions of Theorem A we have 
$f^{-1}(B(f))
\setminus B(f)=\emptyset$ i.e. $B(f)$ is completely 
invariant.) 

 Denote again 
the domain of $\C\setminus \partial$ containing $p$ 
by $D_a$. Denote the complementary open 
topological disc 
$\C\setminus\cl D_a$ by $D'_a$. 

\

{\bf Observations}

3. There is at most one  critical value for $f$ in $D'_a$. So the components of $f^{-1}(D'_a)$ are 
topological discs $D'^j, j=1,...,\hat d$ where $\hat d\le 
d$, with closures in $D'_a$. (In particular $f^{-1}(\cl 
D_a)$ is connected hence $B(f)$ is completely invariant.)

4. Closures of $D'^j$ intersect or "self-intersect"  only at critical points of $f$  and $ \cl  f^{-1}(D'_a)$ is connected , see Fig 1..

\

If the latter were false  then $f^{-1}(D_a)$ would contain 
a nonsimply-connected component $V$ .  But $f$ maps maps $V$  onto the disc $D_a$ so it $V$ would contain a critical point hence $D_a$ would contain a critical value. This would contradict the assumption that all critical values are on the level $a$.

\

After a small perturbation moving (exposing) 
critical values of some $d-1$ critical points towards 
below $a$,  the set $ \cl  f^{-1}(D'_a)$ consists of $\hat d$ closed discs intersecting one another at at most 1 point, which union is connected and simply-connected, Fig 2.

So for $\e>0$ small enough the set  $f^{-1}
(\partial D_{a-\e})$ (where $D_{a-\e}:=D_a\cap\{G<a-
\e\}$), 
is a topological circle $\partial_\e$ bounding a 
topological disc $U\ni p$
and  under $f$ it winds $d$ times onto $\partial D_{a-\e}$. Of course $f(\cl U)=\cl D_{a-\e}\subset U$.

After performing Blaschke type perturbation 
and a holomorphic change of coordinates  on $\C$ 
we arrive at 
a polynomial.\hfill    $\clubsuit$

\

\

\

\

\

\

\

\

\

\

\

\

\

\ \ \ \ \ \ \ \ Figure 1.\ \ \ \ \ \ \ \ \ \ \ \  \ \ \ \ \ \ \ \  \ \ \ \ 
  Figure 2. Critical points 

\ \ \ \ \ \ \ \ \ \ \ \ \ \ \ \ \ \ \ \ \ \ \ \ \ \ \ \ \ \ \ \  \ \ \ \ \ \
\ \ \   $c_1$ and $c_4$ are exposed towards $p$.

\

{\bf Proof of Corollary B.} This Corollary follows from Theorem A because its assertion is an open property in $H^d$ and it is true and easy for polynomials. 

Indeed if $f$ is a polynomial, then with the help 
of a small perturbation we guarantee that 
for each critical value $v_j, j=1,...,d-1$ (different  
from $\infty$) the trajectory $\g_j$ for $\grad G$ 
where $G$ is Green's function in the basin of attraction to 
$\infty$, pole at $\infty$, goes from $v_j$ up to $\infty$. 
In other words $\g_j$ does not go to any critical point 
for $G$. Then there are no critical values for $f$ in 
the topological disc 
$$
U=\C\setminus A \ \ \ \hbox{where}\ \ 
A=\bigcup_{j=1,...,d-1,\  n\ge 0} f^n(\g_j).
$$
Because $f(A)\subset A$, we have a collection of 
branches $g_j:U\to U$ of $f^{-1}$. Denote $g_j(U)$ by 
$U_j$. Then each $z\in J(f)$ is coded by the sequence 
of 
symbols $j_n, n=0,1,...$ where $f^n(z)\in U_{j_n}$.

For each sequence $(j_n)$ the family of maps 
$\Psi_n=g_{j_0}\circ g_{j_1} \circ ... \circ  g_{j_n}$ is 
a normal family of maps on $U$. It is easy to find a 
slightly smaller topological disc $U'\subset U$ so that 
$U'\supset J(f)$ and each 
$$
g_j \ \hbox{ maps} \ \ \cl U' \ \ \hbox{into}\ \  U' .    \leqno (3)
$$
Hence the sequences $\Psi_n|_{U'}$ converge 
uniformly to points in $J(f)$. This proves that the coding is 
one-to-one.
\hfill $\clubsuit$

\

\

{\bf Remark 1.1.} In fact a proof of Corollary B is contained in Proof of Theorem A.  Indeed already Observation 3. gives a partition of the Julia set $J$ into $J\cap D'^j$  and $J\ni z\mapsto (j_n)$ such that $f^n(z)\in D'^{j_n}$ gives a conjugacy to the one-sided shift space of $d$ symbols.

Remark that for this proof of Corollary B, there was no need to refer to Measurable Riemann Mapping Theorem.  Namely there was no need to integrate every time a new conformal structure to obtain a rational mapping. One could work just with smooth maps until the properties described in Observation 3. are  reached.

\

{\bf Remark 1.2.} One can prove Theorem A by perturbing $f$ along a  curve $\hat\g$ close to the Morse curve, dragging the critical value to a small neighbourhood of $p$ just in one step instead of doing it step by step decreasing $s(f(c))$. The property (1) is satisfied by Observations 1. and 2. (Every point of $\hat\g$  can come close to  $\hat\g$ under $f^j$ only further, i.e. if $f^j( \hat\g(\tau))$ is close to $ \hat\g(\tau ')$ then $\tau ' >\tau$. The more so for $f^j$ replaced by 
$h^j_t$ .)

\

{\bf Remark 1.3.} Makienko [M] proved 
the following Proposition which corresponds to our 
Main  Lemma:

If all critical values $v_1,...,v_m$ for critical points in the basin 
$B(f)$ of immediate attraction to an 
attracting fixed point $p_f$  have 
disjoint forward orbits then there exists a 
topological disc $U\subset B(f)$ 
containing $p_f$ such 
that $f|_U$ is injective, $f(\cl U)\subset U$ and 
all $v_j, j=1,...,m$ belong to the annulus $U\setminus 
f(U)$.

\

One can easily change $U$ so that all $v_j$ belong to  
$\p f(U)$. There is a quasiconformal  conjugacy 
$\Phi$ of $f$ on $U$ to $z\mapsto \l z$ on the unit disc 
$\D$, \ $|\l|<1$. So as usually one can pull-back the 
standard conformal structure on $\D$ by $\Phi^{-1}$ 
to $U$ spread it by $\Phi^{-n}, n=1,2,...$ to the whole 
basin of $\p_f$ and complete on the complement of the 
basin  by the standard structure. After integration of this structure the new map satisfies the properties asserted in Main Lemma.

Thus:  Makienko's Proposition (preceded by a small 
perturbation) + Measurable Riemann 
Mapping Theorem, give Main Lemma. 

Conversely, $g$ from the assertion of Main Lemma 
is conjugate to $f$ (provided $f$ is already 
after the first perturbation in Proof of Main Lemma). 
So Makienko's $f(U)$ can be defined as the image 
under the conjugacy of the disc bounded by $\partial$. 
Observe that for that we did not need to construct the 
new holomorphic structures along the proof. (Compare this with Remark 1.1.)

Thus Makienko's Proposition is a topological heart of 
Main Lemma.

\

{\bf Remark 1.4.} Neither in Main Lemma nor in Theorem A 
one needs to assume $f\in H^d$. Proofs work for an 
arbitrary rational map $f, \deg f\ge 2$ if one replaces 
in the statements $H(f)$ by $\Teich_J$. 

The latter 
denotes  the set 
of all such rational maps $g$ that there exists 

\noindent   1) a 
connected open  domain 
${\cal P} \ni 0$ in the complex plane, 

\noindent   2)  a family of quasiconformal homeomorphisms $h_\l:\C\to\C, \l\in \cal P$ with $h_0=\id$
such that for each $z$ the point $h_\l(z)$ depends holomorphically on $z$, 
(such a family is usually called a holomorphic motion ).

\noindent  3)  a family of rational maps $f_\l, \l\in \cal P$ 
with $f_0=f$ and $f_{\l_0}=g$ for some $\l_0\in \cal P$
such that for every $\l$ the map 
$h_\l$ conjugates $f$ with $f_\l$ between $\e$-neighbourhoods of their Julia sets ($\e$ not depending 
on $\l$).

\

\

{\bf Proof of Theorem C.} One proves first an analog to 
Main Lemma:

\

If $f(z)=z+a(z-p)^{t+1}+o((z-p)^{t+1}), \  a\not=0$ then 
$H_1(z)={\l\over z^t}$ for some $\l\not= 0$ conjugates 
$f$ on a "petal" $P\subset B(f)$ to $F(z) = z+1+o(1)$ 
for $z$ with large real part. (See [DH] for the precise 
description.) Next conjugate smoothly 
$F$ to $z\mapsto z+1$ by 
$H_2$.  Define $\Phi_f =H_2\circ 
H_1$ on $P$ and extend it to $B(f)$ by 
$\lim_{n\to\infty} \Phi_f \circ f^n(z) - n$.  Define $G_f(z):=
{1 \over \exp\Re\Phi_f(z)}$ on $B(f)$.

Then after a small perturbation such that 
all critical points in $B$ are simple and  their forward orbits are disjoint, one can find a quasiconformally conjugate $g\in Q^d$ such that all critical values $v_j$ of $g$-critical points in $B(g)$ are in the component  of
$\{G_g=a\}$ bounding a "petal".  (See Fig. 3.)

\

\

A proof of this parabolic version of Main Lemma is the 
same as that of Main Lemma except there are no $M$-
critical points for $G_f$ in $B(f)$. One can think of $M$ 
critical points as belonging to $\p B(f)$, being precisely the set $\p B(f) \cap   \bigcup_{n\ge 0}f^{-n}(\{p\})$.
So the perturbation is not along Morse curves in basins 
$A(q) $  but along curves in $B(f)\cap A(q)$ 
joining directly consecutive 
$S$-critical points in $\p A(q)$.

\

\

Repeat now Proof of Corollary B:

(The assumptions imply $t=1$.)

For $f$ already as $g$ above we take 
 curves $\g_j, j=1,...,2d-2$ joining $v_j$ to $p$. 
We can take as $\g_j$'s the $\Phi^{-1}_f   $-preimages  
of horizontal lines in $\C$. Then $f(\g_j)\subset \g_j$.
So for  $A=\bigcup_{j=1}^{d-1}
\g_j$ the set $\C\setminus A$ is a topological disc 
and we have for the branches of $f^{-1}$ \ 
$g_j(U)\subset U, \ j=1,...,d$.

We find $U'\subset U$ such that (3) holds, except 
$\p U' \cap \p U =\{p\}$ and $g_{j_0}(p)=p$ for the 
branch $g_{j_0}: g_{j_0}(p)=p$, see Figure 3 
 By uniform $(g_{j_0}^n)'
\to 0$ on $\cl U'$ we obtain $\bigcap_n g_{j_0}^n(U')$ 
is only one point $p$. \hfill$\clubsuit$

\

\

\

\

\

\

\

\

\

\

\

\

\

\

\

\

\centerline {Figure 3.}

\

\

{\bf Remark 1.5.}  \  As we mentioned in Introduction the 
fact that under the assumptions of Corollary B or Theorem 
C \   $J(f)$ is a Cantor set is much easier than the fact 
$f|_{J(f)}$ is conjugate to the one-sided shift. 

Indeed, it is easy to prove that $f^N|_{J(f)}$ is conjugate 
to a one-sided shift, for an integer $N>0$. 

Just take a 
small topological disc $D\subset B(f)$ so that $D\ni p$ 
(or a small petal in the parabolic case) such that 
$f(\cl D)\subset D$. , 

For each critical value $v_j$ and every $n\ge 0$ such that $f^n(v_j)\notin D$
join  $f^n(v_j)$ to $D$ by an embedded curve $\g_{j,n}$ 
so that these curves are mutually disjoint. 
Take $N$ such that for each $j,n$ \ \ 
$f^N(\g_{j,n})\subset  D$. Consider now 
$U=\C\setminus A$ where $A=D\cup 
\bigcup_{j,n}\g_{j,n}$. Next proceed as in Proof of Corollary B where $f$ is replaced by 
$f^N$.\footnote{$^4$}{I owe this proof to K. Bara\'nski.}

\

{\bf Remark 1.6.} It is also easy to prove that $f|_{J(f)}$ 
is conjugate to a one-sided topological Markov chain.

Indeed, let $D_0=D$ as in Remark 1.5 . Let $D_n, \ 
n=1,2,...$ be defined recursively:  $D_n$ is
 the component 
of $f^{-1}(D_{n-1})$ containing $D_{n-1}$.

Let $N\ge 0$ be large enough that there are no 
critical values in $\C\setminus D_N$.

Let $U_1, ...,U_K$ denote all the components of 
$\C\setminus \cl D_N$. These are topological discs 
because $D_N$ is connected. Consider the family 
of topological discs 
$g_j(U_k)$ for all branches $g_j, j=1,...,d$ of $f^{-1}$ 
and $k=1,...,K$. Then  for every sequence of pairs 
$(j_n,k_n), n=0,1,...$ 
such that $g_{j_n}(U_{k_n}) \subset U_{k_{n+1}}$ 
there exists precisely one point  $z\in J(f)$ such that 
for every $n\ge 0$ \ $f^n(z)\in g_{j_n}(U_{k_n})$.

\

{\bf Remark 1.7.}\footnote{$^5$}
{Before reading this  Remark the reader is 
advised  to read Section 2.} \  
Considering the situation in Theorem A
such that $p$ is already not critical (say all $f$-critical points belong to $B(f)$ and $\Arg 2\pi f'(p)$ is rational), John Milnor asked 
whether it is possible to find a set $A$ which is union 
of parts of trajectories of $\grad G$  contains all 
critical values is compact in $B(f)$ is $f$-forward invariant
and is connected simply-connected (a tree).

This would immediately allow to prove Corollary B with 
the set  $A$ as here.

Unfortunately the answer is negative. We can modify 
an exotic example in Section 2  
(the case $d=3$, Fig. 5) so that $f(c_4)=f^2(c_3)=a$  
the pole. Then the set of saddles for $G$ is 
$S=\bigcup_{n\ge 0}f^{-n}(c_2)$. So every curve 
$\Gamma $ built from the pieces of trajectories of 
$\grad G$ joining $c_4=f(c_4)$ to $\infty$ passes 
a point of $f^{-n}(\{c_2\})$. Hence $f^n(\Gamma) $ 
joins $\infty$ with $\infty$ passing through $c_2$, 
hence it is a loop, i.e. $A$ is not a tree.

(Even the assumption $f(c_4)\not=a\not=f^2(c_3)$ 
does not help if  $f(c_4), f^2(c_3)$ are close to $a$. 
$\Gamma$ must still leave the basin $A(a)$ for 
$\grad G$ passing through a point in 
$\bigcup_{n\ge 0}f^{-n}(\{c_2\})$.)

\

{\bf Remark 1.8.}$^5$ \   Though in Main Lemma we can arrange 
all critical  values of critical points in $B(f)$ on one 
component of a level of $G$ it sometimes is not so 
for critical points (unless all critical points for $f$ 
or all but one, are in $B(f)$ as in Theorem A, see Fig 2.). 
Again we modify an exotic example from Section 2. 
Here degree of $f$ is 5. Start with a cubic polynomial 
$P$ which has degree 1 on $\p A_1$,  degree 2 
on $\p A_2$, see Fig 3, Section 2, and $P$ maps 
the critical point $a\in A_2$ to the critical point at 
the self-intersection of the figure 8 line $\p A_1\cup \p 
A_2$, which escapes to $\infty$. 

Consider now the function $z\mapsto P(z)
+{b\over (z-a)^2}$ with $b$ small real positive number.
We obtain the picture as on Figure 4. Do the surgery 
as in Section 2 to have it holomorphic. For a final example split $a$ into two different poles which gives an 
$f$-critical point $c_8$ between them, Fig. 4. Make 
Blaschke type perturbation close to $\infty$ to have $\infty$ not critical  and move $f(c_8), f(c_1), f(c_2)$ to one level with $f(c_3)$.

\

\

\

\

\

\

\

\

\

\

\

\

\centerline{Figure 4. $f^2(c_4)=c_4, f^2(c_6)=c_6, f(c_5)=c_5, f(c_7)=c_7$.}

\

\

{\bf Section 2. \  Exotic basins.} 

\

{\bf Proof of Theorem D,  case $d=3$.} 
We start with a geometric description of an exotic example of degree 3, illustrated on Fig. 5 .

\

Start with a quadratic polynomial $P(z)=z^2+c$ with the  critical point $c_2=0$ , escaping to $c_1=\infty$ a attracting fixed point 
of  multiplicity 2. 
The level $\partial=\{G=t\}$ 
of Green's function of the basin of attraction 
to $\infty$ with the pole at $\infty$, containing $c_2$, 
is  figure eight.  Now we change the map 
on  $B_2$, one of the two discs $B_1, B_2$ 
bounded by $\partial$ as follows:
  
Draw two little discs $D_1, D_2$ in $B_2$, intersecting one another. Let $D_1\cap D_2$ be maped 1-to-1 onto $\C\setminus (B_1\cup  B_2)$ . Let $D_1\setminus D_2$ goes
onto $B_1$ and  $D_2\setminus D_1$ goes
onto $B_2$ both both proper maps 
with degree 2. So there are critical points $c_3\in D_1\setminus D_2$ and 
$c_4\in D_2\setminus D_1$. On $D_2$  this map $f$ is quadratic-like so we 
can do  anything  there, for example $f(c_4)=c_4$. 
On $D_1$ the map $f^2$ is quadratic-like so we can assume $f^2(c_3)=c_3$.

\

\

\

\

\

\

\

\

\

\

\

\

\

\

\centerline{Figure 5.}

\

The rational function is obtained out of this topological picture by the quasi-conformal surgery technique [D].
We shall explain it closer now:

\

We need following Lemma which generalizes 
Douady-Hubbard's theorem that a po\-ly\-no\-mial\--like 
mapping is quasiconformally conjugate to a polynomial 
[DH1]:

\

{\bf Lemma 2.1.} \  Let $U\subset \C$ be an open set 
(not 
necessarily connected or simply-connected) with 
boundary being a family of smooth Jordan curves.
Let $F_1:U\to U$ be a holomorphic map such that its  $F_1(\cl U)\subset U$. 
(We denote the continuous  extension of $F_1$ to 
$\cl U$ by the same symbol  $F_1$.)  

Let $V\subset \C$ be homeomorphic to $\C\setminus \cl 
U$ by a  homeomorphism $h_1$ which extends 
orientation 
preserving to a homeomorphism of $\C$. (Again we do not assume $V$ is connected or simply-connected.) 
Let $F_2:V\to 
\C$ be a holomorphic map.  We also suppose that 
the boundary of $V$ is smooth and denote by $F_2$ the continuous extension of th original $F_2$ to $\cl V$.

Suppose the family of curves being the components 
of $\partial U$ has the same combinatorics in $\C$ 
as the family 
of curves being the components of $\partial V$. 
We mean by this, that

1.  There exists a homeomorphism $h_2:  \p U \cup \p 
F_1(U) \to \p V \cup \p F_2(V)$  such that the boundary 
of each component of $\C\setminus (\p U \cup \p 
F_1(U)$ is mapped to the boundary of a component of 
$ \p V \cup \p F_2(V)$

2.  For each component $\p$ of $\p U$ the map $h_2$ 
maps $\p$ to $h_1(\p)$, \  $F_1(\p) $ to $F_1(h_1(\p))$ 
and there
a continuous map (a lift) $\tilde h_2:\p \to h_1(\p)$ such that 
on $\p$ we have $h_2 \circ F_1 = F_2 \circ \tilde h_2$.
(i. e. $h_2$ preserves orders between $F_1(\p)$ and $F_2(h_1(\p))$).

\

Then there exists a rational map $f:\C\to\C$ and an 
open $W\subset \C$ such that $f$ is quasiconformally conjugate to $F_ 1$ on $W$ and quasiconformally conjugate to $F_2$ on $\C\setminus W$.

\

{\bf Proof.} We  replace $h_2$ on $\p U$ by the lift
$\tilde h_2$ and then extend 
it from $\p U \cup \p 
F_1(U)$ to a quasiconformal homeomorphism $h:\C\to 
\C$. Define $F:\C\to\C $ by $F_2$ on $\cl V$ and by 
$h\circ F_1 \circ h^{-1}$ on $\C\setminus V$.

 Let $\mu_0$ denote the standard conformal 
structure on $\C$. Take  $\mu_1=h_*(\mu_0|_U)$ on 
$\C\setminus \cl V$. (Think about $\mu_1$ as a field of 
ellipses, up to a multiplication by a positive function.)
For each $z\in V$ define $\mu_1(z)$ as a pull-back 
  $F^{-n}_*((\mu_1)(F^n(z))$ where $n\ge 0$ is such 
that $F^n(z)\in \C\setminus \cl V$. If such $n$ does not 
exist take $\mu_1(z)=\mu_0(z)$. This is correct due to the crucial property $F(\C\setminus \cl V)\subset \C\setminus \cl V$. As $F_2$ is holomorphic, $\mu_1$ is in $L^\infty$ !

Now integrate $\mu_1$. In the new coordinates 
$F$ changes to a rational map $f$ we looked for. \hfill 
$\clubsuit$

\

\

Now we construct $F_1$ and $F_2$ satisfying the assumptions of Lemma. It is illustrated on Fig. 6.

Take the polynomial $P(z)=z^2+c, \ c<-2$. Make 
$F_1$ by adding to $P$ a term ${b \over z-a} $ for 
$a=\sqrt{-c}\in P^{-1}(0)\cap  
B_2$, \ $(a>0)$. Let $b>0$ be 
small so that for our $F_1$ the level $\hat\p =\{\hat G = 
t_0\}$ 
containing the $F_1$ critical point $\hat c_2$ close to 
$c_2=0$ is figure 8 close to $\p$.  
Here $\hat G$ is defined 
analogously to Green's function or to $G$ in Section 1:  
on the basin of 
attraction  to $\infty$ by $F_1$, one  defines 
$\hat G(z)=\lim_{n\to\infty}2^{-n}\log |F_1^n(z)|$. 
Denote discs bounded by adequate parts of $\hat\p$ 
close to $B_1, B_2$ by $\hat B_1, \hat B_2$ respectively.

It is easy to compute that $\hat c_2={b\over 
2a^2}+o(b)$ and for two other 
$F_1$-critical points $\hat c_3, \hat 
c_4$ we have  $F_1(\hat c_{3,4})=
\mp 2\sqrt{2a}\sqrt{b} +o(\sqrt{b})$.
So $F_1(\hat c_3)\in \hat B_1, F_1(\hat c_4)\in \hat B_2$.

Let $0<t_2<t_1<t_0$ with $t_2\approx t_1\approx t_0$
and denote by $K_2, K_1$ the topological discs both in 
$\hat B_1$ bounded by $\{\hat G=t_2\}$, respect. $\{\hat G=t_1\}$. Denote by $K_4'$ the topological disc in 
$\hat B_2$ bounded by $\{\hat G=t_1\}$. Finally denote 
the component of $F_1^{-1}(K_1)$ in $B_2$ by $K_3$ 
and denote the component of $F_1^{-1}(K_4')$ in $B_2$
by $K_4$. 
We have $\hat c_3\in K_3, \ \hat c_4\in K_4$.

Define $U:=\C\setminus (K_2\cup K_3 \cup K_4)$.
We have $F_1(\cl U) \subset U$. So $F_1$ and $U$ 
satisfy the assumptions of Lemma. Now we need to 
define $F_2$: 

Set  $F_2(z):=z^2$ on a  geometric 
disc $L_4=\{|z|<r_4\}, \ r_4>1$. 
Take a disc $L_3= \{|z-z_0|<r_3\} \subset F_2(L_4)
\setminus \cl L_4$ 
and define 
$\tilde F_2(z)=(z-z_0)^2+z_0$. One finds large 
$r_4, r_3$  such that 
$\tilde F_2(L_3)\supset \cl F_2(L_4)$. 
Pick in $\tilde F_2(L_3)\setminus F_2(L_4)$ two 
discs 
$L_2\subset L_1$ of the form 
$L_1=\{|z-z_1|<r_1\}, \ L_2=\{|z-z_1|<r_2\}, 
r_2<r_1$.
Take an affine holomorphic map $\Psi:L_1\to 
\tilde F_2(L_3)$ (onto). Define 

$F_2=\Psi^{-1}  \circ  \tilde F_2$ on $L_3$ and 

$F_2=\Psi |_{L_2}$ on $L_2$.

\noindent   We care to have $r_2$ so close to $r_1$ that $\Psi(L_2) 
\supset \cl F_2(L_4)$. 

Now take $V=L_2\cup L_3 \cup L_4$ and $F_2$ 
defined on   $V$ as above.

The assumptions of Lemma are satisfied.  So we can 
"glue"  $F_1$ and $F_2$ in one rational mapping $f$.\

\

\

\

\

\

\

\

\

\

\

\centerline {Figure 6.}

\

\

Observe finally that $J(f)$ is disconnected because 
$F_1^n(\hat c_2) \to \infty$  and moreover 
$F_1^n(l) \to \infty$ where $l=\{\Im z = \Im \hat c_2 \}$. 
The line $l$ separates $K_2$ from say $K_4$.  
Both $K_1$ and $K_4$ intersect $J(f)$  (in the 
coordinates after the integration of $\mu_1$) so 
the intersections belong to different components of 
$J(f)$. 

The degree of $f$ on the basin of attraction to $\infty$ is 
3 because such is the degree of $F_1$ on $U$. Only 
two critical points: $\infty$ and that one corresponding to 
$\hat c_2$ belong to the basin, because $\hat c_ {3,4} $ 
do not escape under the iteration by $F_2$. Theorem D 
is proved for $d=3$.

\

{\bf Remark 2.2.} Observe that in appropriate holomorphic coordinates on $\C$ we have $f(z)=z^2+c+
{b \over z-a}$. Indeed after subtracting from $f$ 
constructed above the principal  part of the 
Laurent series expansion at the pole, we are left with 
a quadratic polynomial. By an affine holomorphic change of coordinates we arrive with teh polynomiial to $z^2+c$.

\ 

{\bf Remark 2.3.} One should be careful in the above construction because not every branched cover 
of $\C$ preserves a conformal structure.  
Above, an annulus ${\cal A}$ in $B_2$ containing $c_3$ and 
$c_4$ is mapped in a proper way by $f$ to the disc 
$D'=\{G<t'\}, t'>t$ ,  i.e. a disc containing $\infty$ , 
outside  the figure 8 level $\{G=t\}$. 

Instead of mapping $c_3$ into $D_1$ so that 
$F^2(c_3)=c_3$ we can map $\cal A$ onto $D'$ in 
a proper way so that $f(c_3)=c_3$ and $f(c_4)=c_4$. 
This will be a topological branched cover.  However 
it does not allow a holomorphic invariant structure. 

If it allowed, for $c_3$  close to $c_4$, 
 for ${\cal A}$ small but of a definite modulus, 
then in the limit 
after rescallings of  ${\cal A}$'s to be of a definite size 
we would end up with a covering map of an annulus 
with two punctures to a disc with a puncture (covering 
without branching point). This is not possible by the Euler 
characteristics argument. It means $c_3$ and $c_4$ 
cannot  be too close in $D'$. 
 
Another argument is that such $f$ would have 3 
superattracting  fixed points. So it would be a Newton's 
method rational function of a degree 3 polynomial, see 
Introduction. 
But the basin of attraction to  $\infty$ is not simply-connected. This contradicts a theorem that the basins of immediate attraction to the
attracting fixed points for Newton's method are simply-connected [P].

\

\

{\bf Proof of Theorem D, the general case $d\ge 3$.}

\

We shall realize holomorphically the picture on Fig. 7:

\

\

\

\

\

\

\

\

\

\

\

\

 \centerline{Figure 7.}

\

\

On Fig. 7, \ 
 $D_j$ is mapped properly on $B_j$ for $j=1,2$. 
Each $D_j$ contains $d-2$ critical points. The points 
$a_1,...,a_{d-2}$ are poles.

\

We proceed similarly as in the case $d=3$.  Let  
$$
F_1(z)=z^2 +c + b(\sum_{m=1}^{d-2}{1  \over z-a_m}).
$$

Take $a_m=\sqrt{-c} + imT$ for a real constant $T: 
0<T\ll 1$ in particular $T$ small enough that all $a_m$ 
are well in $B_2$.

\

For $b$ real $b>0$ small, there is a small annulus around 
each pole $a_m$, containing two critical points 
$$
\hat c_{m,3}, \hat c_{m,4}= a\mp\sqrt{b/2a} + 
o(\sqrt{b}).$$

The corresponding critical values 
$$
v_{m,3}=F_1(\hat c_{m,3}), \ \ 
v_{m,4}=F_1(\hat c_{m,4})
$$
are
$$
a_m^2+c\mp 2\sqrt{2a}\sqrt{b} +o(\sqrt{b}) =
2imT + (-m^2T^2 \mp 2\sqrt{2a}\sqrt{b} +o(\sqrt{b}) ).
$$

(Computing $ \hat c_{m,3(4)}$ and $v_{m,3(4)}$ it is 
comfortable to consider $z^2+c+{b\over z-a_m}$.
Other terms ${b \over z-a_t}$ have only the $O(b)$ influence.)

\

The $F_1$ critical point close to 0 is $O(b)$. 

\

Take $K_1, K_4'$ from the case $d=3$ slightly modified, 
larger than the original ones:  
Let $l$ be the line (parabola) $(2iT\tau,  -T^2\tau^2)$ 
for   
$\tau>0 $. Observe that the critical values $v_{m,3}$ 
are to the left of $l$, and $V_{m,4}$ to the right of $l$. 
We extend $K_1, K_4'$ to $\hat K_1, \hat K_4'$ almost 
to $l$  to capture $v_{m,3}, v_{m,4}$ respectively.
See Figure 8.

\

\

\

\

\

\

\

\

\

\

\

\

\

\

\centerline{Figure 8.}

\

\

Consider the topological discs 

$K_{m,3}=\Comp F_1^{-1}(K_1), \ \  K_{m,4}=\Comp 
F_1^{-1}(K_4')$, where $\Comp$ means the component 
containing $\hat c_{m,3}, \ \hat c_{m,4}$ respectively.

$K_2$ is as in the case $d=3$, such that 
$F_1(K_2)\supset \cl K_4' .$

Finally set

$$ U:= \C\setminus (K_2\cup \bigcup_{m=1}^{d-2} 
(K_{m,3}\cup K_{m,4}).$$

\

The rest of the construction of $f$ is the same as 
for $d=3$. When we make quadratic-like maps 
$\tilde F_2$ on $L_{m,3}$ and $F_2$ on $L_{m,4}$ 
we have a complete freedom of which quadratic polynomials we glue in, in particular of whether we 
want the corresponding critical points to escape or not. 
(In particular if no critical point escape we have the most surprising case $k=2$ of the assertion of Theorem D.)
\hfill$\clubsuit$

\

For the completeness of the exposition we shall prove 
the following simple facts (The first of them stated already 
in Introduction) :

\

{\bf Proposition 2.4.} a) \  Let $f\in H^d$ be a polynomial 
with $B(f)$ the basin of attraction to $\infty$ not simply-connected. Then $B(f)$ contains at least one critical 
point different from $\infty$.

b) \   More generally, if $f\in Q^d$ and for a pair of 
topological discs $A, A_1: \ \cl A \subset A_1$ the map 
$f|_A:A\to A_1$ is proper of  degree $d' \le d$, then 
$A_1$ is in the basin of attraction $B(f)$ to an 
attracting fixed point and if $B(f)$ is not simply-connected 
then it contains at least $d'$ critical points. This concerns in particular the case $d'=d$ in which $f|_{\C\setminus\cl 
A_1}:\C\setminus\cl A_1 \to \C\setminus\cl A$ is polynomial-like.

\

{\bf Proof.} \  a) \  Take a topological disc $D=\{G > 
a\}$  around $\infty$ (cf. Proof of Corollary B or Remarks 
1.5, 1.6). If there are no critical points in $B(f)$ (except $\infty$) then $f^{-n}(D)$ is an increasing sequence of topological discs, so $B(f)=\bigcup_{n\ge 0} f^{-n}(D)$ is a topological disc, hence $B(f)$ is simply-connected. 

(Remark that we already used the argument, that if there 
is only one critical value for a proper map $f:W_1\to 
W_2$ where $W_2$ is a topological disc , then $W_1$ 
is also a topological disc, in Remark 2.2)

b) \ The proof is similar.  There are $d'-1$ critical points in $A_1$ and there must be a critical point in $B(f) \setminus A_1$. 

(One can also deduce b) from a) using Blaschke type 
perturbation, Section 1.)
\hfill$\clubsuit$

\

{\bf Proposition 2.5.} Every non simply-connected 
immediate basin 
of attraction to an attracting or parabolic fixed point 
(with $f'(p)=1$) contains at least 2 different critical values of critical points 
in the basin. 

\

This complements
 Theorem D: the integer $k$ cannot be 
less than 2.

\

{\bf Proof.}  Consider the sets $D_n$ defined in Remark 
1.6.  As $B(f)=\bigcup_{n\ge 0} D_n$ is not 
simply-connected, there exists $n$ such that $D_n$ is simply-connected and $D_{n+1}$ is not. Then $D_n$ 
contains at least two different critical values of critical 
points in $D_{n+1}$.\hfill$\clubsuit$

\

(Almost the same argument proves the above for periodic basins, period larger than 1.)

\

{\bf Proof of Corollary E.} Let $f\in H^d$  be as in Theorem D for $d\ge 3, k=2$. Consider an 
arbitrary  $g\in H(f)$. Then there exist a real  continuous 
1-parameter family of homeomorphisms $h_t:\C\to\C$  
and a real 1-parameter family of maps $f_t\in H(f)$ 
having the same properties as $h_\l$ and $f_\l$ in 
Remark 1.4. (There exist complex families precisely 
as in Remark 1.4, but we do not need them here.)
Then $\sharp(\Crit (f_t)\cap B(f_t)) $ is constant because 
critical points cannot be too close to $J(t)$, where $|f'|>0$ 
uniformly, they move continuously with $t$, so they cannot jump between components of $\C\setminus J(f_t)$. So by Proposition 2.4, $B(g)$ cannot be the basin of attraction to $\infty$ for a polynomial. But degree of 
$f$ hence $g$ on each other invariant 
basin $B_1$ is less than $d$. (Otherwise $\p B_1=J(f)$ would be connected and it is not because $B(f)$ is not simply-connected.)
So 
$g$ cannot be a polynomial.\hfill$\clubsuit$

\

Remark that it follows  from Proposition 2.4,b) and 
above Proof that none $g\in H(f)$ 
has a polynomial-like restriction of degree 
$d$.

\

\

\

{\bf Section 3. A 1-parameter family of functions 
joining an exotic $z\mapsto z^2+c+{b\over z-a}$ to 
Newton's method rational function.}

\

Let $f(z)=z^2+c+{b\over  z-a}$. Then $f'(z)=2z-{b\over 
(z-a)^2}$. The equation for the critical points in $\C$ is 
$$
2z(z-a)^2=b
$$
Suppose that $w=c_4$ is an $f$-fixed critical point, see 
Fig 5, Section 2.
(This restricts the number of parameters to 2.)
We obtain
$$
\ \ \ \ \  (w^2-w+c)(w-a) =-b
$$
$$
\ \ \ \ \  2w(w-a)^2=b
$$

Let $a=kw$. We parametrize $f$ by $k$ and $w$. We obtain:
$$a=kw, \ \ b=2w^3(1-k)^2, \ \  c=w^2(2k-3)+w$$

The critical points are $u=c_2, v=c_3, w$, where 
$$
u,v=w(-{1\over 2} +k  \mp{1\over 2}\sqrt{4k-3})
$$

\

Given a parameter $k<1$ sufficiently close to 1, one finds 
$w$ such that $f^2(v)=v$ and the trajectory of $u$ escapes to $\infty$. 
For  $k= .85$ one finds $w\approx 1.88053$.
This is an exotic example as in Figure 5, Sec.2. 
The reason is that the geometry is as in Fig. 6, Sec.2, so the basin $B$ of attraction to 
$\infty$ is connected, i.e. the immediate basin is completely invariant. 
The picture is similar to that in Fig 10d.

(It is not clear to me whether just the escape of $u$ to $\infty$ proves the connectedness of $B$. One should be careful because for $f$ being only a topological branched cover this is 
not so, see the example in Remark 2.3.) 

\

In the rest of this Section we discuss the change of dynamics for 
varying parameter $w$.

 \

For $k=.85$ the number $w=1.88053$ is in the principal part of a 
Mandelbrot-like set $M(c_3)$, symmetric with respect 
to real $w$'s, pronged to the left. For $w\in M(c_3)$, Julia set for quadratic-like $f^2|_{D_1}$ (see Fig. 5, 
Section 2)  is connected and we still have exotic maps. 

 Now let us decrease $w$. It 
leaves $M(c_3)$ at $w\approx  1.86874$ and below that $w$ the trajectory of  $v$ escapes from $D_1$.
It need not escape to $\infty$. There is a sequence of 
intervals where $f^{2n}(v)$ hits $B_w$ the basin of immediate attraction to $w$, \ \ $n$ decreases to 2. 
Later on, after escape, again $f^4(v)\in B_w$ but 
$f^2(v) <u$ (before, it was between $w$ and $v$). 
This happens  at 
$w\approx 1.63045$. See Fig. 9.

\

\

\

\

\

\

\

\

\

\

\

\

\

\

\

\centerline{Figure 9. $k=.85, w=1.63045$.}

\

\

At some parameter $w$ the trajectory of 
 $u=c_2$ stops to 
escape to $\infty$. It hits $B_w$. But next with further decrease of $w$ it can again escape to 
$\infty$. 

Starting from $w\approx 1.541549$  the trajectory of $u$ 
neither escapes to $\infty$ nor to $w$.
The parameter $w$ is in a Mandelbrot-like set $M(c_2)$ 
pronged towards right.  In fact at this parameter 
$f^2(u)\in (v,a)$. Only after some further decrease of $w$ 
we arrive at $f^2(u)\in (u,v)$, so that one has a unimodal map $f:(f(u),v)\to (f(u),v)$. 

$w\approx .7136114$ is 
in the principal part of $M(c_2)$ and $f$ is Newton's. 
Then $f^n(c_3)\to c_2=f(c_2)$. The number $w\approx  
.301$ is still in the principal part of $M(c_2)$ and $f$ is 
Newton's but now $f^n(c_2)\to c_3=c_3$.\footnote{$^6$}{This 
description comes out of a computer made picture in 9 colors, showing whether $c_2,c_3$ escape to $\infty$ to 
$w$ or make something else.}

\

Let us present now pictures from this experiment for $k=.81$.

On Fig.10, $k= .81$, white is the basin of attraction to $w$, 
grey the basin of $\infty$, black is the complement. 
For Newton's, Fig 10a, black contains both $c_2$ and $c_3$, so it has a connected interior and accesses the only repelling 
fixed point in two channels. Let $w$ grow.    For $w\approx 1.37$  black 
Newton's basin has bifurcated to period 4, Fig 10b. 

For $w\approx 1.4961$, Fig 10c., \ $w$ is already in $M(c_3)$ but $u=c_2$ does not escape to $\infty$. It is in the basin of $w$. The basin of $\infty$ is not connected. This is so because the immediate basin (and the whole 
basin too) contains only 1 critical point: $\infty$. 
So it is simply-connected, see Prop. 2.5. Hence $f$ has only degree 2 on this immediate basin.

For $w\approx 1.51545$ \ $u$ 
escapes to $\infty$. The basin of $\infty$ becomes 
connected. This is one of our exotic examples: see 
Fig. 10d. 

\

\

\

\

\

\

\

\

\

\

\

\

\

\

\

\

\

\

\

\

\centerline{Figure 10a.
$k=.81, w=.63$, window $-2-2i,2+2i$.}

Iteration of Newton's method rational map for a polynomial. Black, white and 
grey are basins of attraction to the zeros of the polynomial. 

\

\

\
\

\

\

\

\

\

\

\

\

\

\

\

\

\centerline{Figure 10b. 
$k=.81,w=1.37$, window $-2-2i,2+2i$.}

Black Newton's basin has bifurcated to period 4 immediate basin and its 
pre-images.

\

\

\

\

\

\

\

\

\

\

\

\

\

\

\

\

\

\

\

\

\

\

\centerline{Figure 10c. $k=.81, w=1.49$, window $-2-2i,2+2i$.}  

The map is still not exotic because the trajectory of the critical point $u$ 
is attracted to $w$.

\

\

\

\

\

\

\

\

\

\

\

\

\

\

\

\

\

\

\

\

\

\centerline{Figure 10d. $ k=.81, w=1.51545$, window $-i,2+i$.}
 
This is an exotic map. The pattern is as in Figures 5,6. The union of 
black and white 
does not separate plane anymore, $u$ escapes to $\infty$.

\

\

\

{\bf Question 3.1.} \  In the set of Newton's method rational functions $NP_\l$ 
for the polynomials $P_\l(z)= z^3+(\l -1)z - \l$ 
there exist Mandelbrot-like sets where the critical point 
different from the zeros of $P_\l$ converges to a periodic  
attracting orbit different from these zeros, [CGS]. Do 
 these sets move to  $M(c_3)$ sections of 
the set exotic maps when we change parameters from 
Newton's to the exotic ones?

\

{\bf Question 3.2.} \ Describe precisely how does the dynamics bifurcate (what is the limit
behaviour 
of the trajectories of $c_2$ and $c_3$)  for real parameters $k,w$. This is
the question on the iteration of the real map having 2 critical point, namely
our $f$ restricted to $(-\infty,a)$. (The right branch from $a$ to $\infty$ 
does not take part in the recurrence because for $z>a$ for every $n\ge 0$ \ 
$f^n(z)\ge a$.)

\

\

\

\

{\bf References }

\

[AB] Ahlfors, Bers, Riemann's mapping theorem for variable metrics, 
Ann. of Math. 72 (1960), 385-404.

[B] B. Bojarski, Generalized solutions of systems of differential equations 
of first order and elliptic type with discontinuous coefficients, Mat. Sb. 43
(85) (1957), 451-503. (In Russian.)

[Boyle] M. Boyle, a letter.

[CGS] J. Curry, L. Garnett, D. Sullivan, On the iteration of a rational 
function: computer experiments with Newton's method, Comm. Math. Phys. 91 
(1983), 267-277.

[D] A. Douady, Chirurgie sur les applications holomorphes, Proc. of the ICM
Berkeley 1986, 724-738.

[DH1] A. Douady, J. Hubbard, On the dynamics of 
polynomial-like mappings, Ann. Sc. \'Ecole Norm. Sup. 
18 (1985), 287-243.

[DH2] A. Douady, J. Hubbard,  \'Etude dynamique des 
polyn\^omes complexes, Publications Mathematiques 
d'Orsay: 2 (1984), 4 (1985). 

[GK] L. Goldberg, L. Keen, The mapping class group of a generic quadratic 
rational map and automorphisms of the 2-shift, Invent. Math. 101 (1990), 335-372.

[M] P. Makienko, Pinching and plumbing deformations of 
quadratic rational maps, preprint  International Centre for Theor. Phys., Miramare -Trieste, February 1993.

[P] F. Przytycki , Remarks on simple-connectedness
of basins of sinks for iterations of rational maps, 
Banach Center Publ. 23 (1989), 229-235. 

\

\

Feliks Przytycki,

Institute of Mathematics, Polish Academy of Sciences.

ul. \'Sniadeckich 8,  00-950 Warsaw, Poland

\end